\documentclass[12pt,reqno]{amsart}
\usepackage{amsfonts,amsmath,amscd,amssymb,amsthm}
\usepackage{enumerate,graphicx,hyperref,perpage,url}
\usepackage[margin=2.5cm]{geometry}

\theoremstyle{plain}
\newtheorem{thm}{Theorem}

\newtheorem{cor}[thm]{Corollary}

\newtheorem{pblm}[thm]{Problem}

\theoremstyle{remark}

\newtheorem*{rmk}{\textbf{Remark}}

\numberwithin{equation}{section}

\MakePerPage{footnote} \allowdisplaybreaks 

\newcommand\Area{\mathrm{Area}\,}
\newcommand{\B}{\mathbb{B}}
\newcommand{\C}{\mathbb C}

\newcommand{\HH}{\mathcal H}

\newcommand{\M}{\mathbb M}

\newcommand{\NN}{\mathcal N}

\newcommand\R{\mathbb R}
\newcommand{\s}{\mathbb S}

\newcommand{\T}{\mathbb T}

\newcommand{\Vol}{\mathrm{Vol}}

\title{Nodal lengths of eigenfunctions in the disc}

\author{Xiaolong Han}
\email{Xiaolong.Han@csun.edu}
\author{Michael Murray}
\email{Michael.Murray.921@my.csun.edu}
\author{Chuong Tran}
\email{Chuong.Tran.561@my.csun.edu}
\address{Department of Mathematics, California State University, Northridge, CA 91325, USA}

\thanks{} 
\date{}

\begin{document}
\maketitle

\begin{abstract}
In this paper, we derive the sharp lower and upper bounds of nodal lengths of Laplacian eigenfunctions in the disc. Furthermore, we observe a geometric property of the eigenfunctions whose nodal curves maximize the nodal length.
\end{abstract}

\section{Introduction}
In an $n$-$\dim$ smooth and compact Riemannian manifold $(\M,g)$, let $\Delta=\Delta_g$ be the Laplacian and $u$ be an eigenfunction with eigenvalue $\lambda$, i.e. $-\Delta u=\lambda u$. If $\M$ has smooth boundary, we impose Dirchlet or Neumann boundary condition. Yau \cite{Y} conjectured that
\begin{equation}\label{eq:Yau}
c_1\sqrt\lambda\le\HH^{n-1}(\NN(u))\le c_2\sqrt\lambda
\end{equation}
for some constants $0<c_1,c_2<\infty$ depending on $(\M,g)$ and independent of the eigenvalues $\lambda\to\infty$. Here, $\HH^{n-1}$ denotes the $(n-1)$-$\dim$ Hausdorff measure and $\NN(u)=\{x\in\M:u(x)=0\}$ denotes the nodal set of the function $u$.

If the metric $g$ is analytic, then \eqref{eq:Yau} was proved by Donnelly-Fefferman \cite{DF1} (see also Lin \cite{Li} for the upper bound); if $g$ is smooth, then Logunov \cite{Lo1, Lo2} recently showed that 
$$c_1\sqrt\lambda\le\HH^{n-1}(\NN(u))\le c_2\lambda^\alpha,$$
in which $\alpha>1/2$ depends only on $n=\dim\M$. There are partial results in this direction \cite{Bru, ChMu, CoMi, DF2, DF3, HL, HaSi, HeSo, HW, SZ1, SZ2} etc, c.f. the survey by Zelditch \cite{Z2}.

In this paper, we are concerned with the precise and sharp dependence of the constants $c_1$ and $c_2$ in \eqref{eq:Yau} on the geometry $(\M,g)$. That is, write
\begin{equation}\label{eq:infsup}
H_1(\M)=\liminf_{\lambda\to\infty}\frac{\HH^{n-1}(\NN(u))}{\sqrt{\lambda}}\quad\text{and}\quad H_2(\M)=\limsup_{\lambda\to\infty}\frac{\HH^{n-1}(\NN(u))}{\sqrt{\lambda}}.
\end{equation}
Thanks to \cite{DF1, Lo2}, we know that $0<H_1\le H_2<\infty$ on analytic manifolds and $H_1>0$ on smooth manifolds. 

Regarding the two limits in \eqref{eq:infsup}, we also pursue the categorization of the eigenfunctions (in relation to the geometry of the manifold) which saturate $\liminf$ and $\limsup$ of $\HH^{n-1}(\NN(u))/\sqrt{\lambda}$ as $\lambda\to\infty$. That is, what geometric properties do the nodal sets of these eigenfunctions achieving the limits in \eqref{eq:infsup} have? Our primary interest is to categorize the sequence of eigenfunctions whose nodal curves are geodesics in the manifold, if these eigenfunctions exist. See Corollary \ref{cor:geod} and Problem \ref{pblm}.

We begin from the easiest case, an interval $\mathbb I_L=[0,L]$. The $j$-th Dirichlet eigenfunction is
$$u_j(x)=\sin\left(\frac{j\pi}{L}x\right)\quad\text{ with eigenvalue }\lambda_j=\left(\frac{j\pi}{L}\right)^2,\ j=1,2,3,...$$
The nodal set (in the interior of the domain) of $u_j$ is a collection of nodal points
$$\NN(u_j)=\left\{\frac Lj\cdot l:l=1,...,j-1\right\}.$$
Hence, the size of the nodal set of the $j$-th Dirichlet eigenfunction $u_j$ on $[0,L]$, i.e. the number of nodal points, is 
$$\HH^0(\NN(u_j))=j-1=\frac{L\sqrt{\lambda_j}}{\pi}-1,\quad j=1,2,3,...$$
Therefore, in \eqref{eq:infsup} we have that
$$H_1(\mathbb I_L)=H_2(\mathbb I_L)=\lim_{\lambda\to\infty}\frac{\HH^0(\NN(u))}{\sqrt\lambda}=\frac1\pi L.$$
One can similarly show the same results for Neumann eigenfunctions in $[0,L]$.

If $\dim\M\ge2$, then very little is known about \eqref{eq:infsup}. To the authors' knowledge, the first result in this direction is due to Br\"uning-Gromes \cite[Equation (8)]{BG}: They remarked that in an irrational rectangle $R$ with side-lengths $a$ and $b$ for which $a^2/b^2$ is irrational, 
\begin{equation}\label{eq:BG}
H_1(R)=\frac1\pi\Area(R)\quad\text{and}\quad H_2(R)=\frac{\sqrt2}{\pi}\Area(R).
\end{equation}
However, $H_1$ and $H_2$ are not known in more general rectangles. In \S\ref{sec:rectangles}, we discuss finding $H_1$ and $H_2$ in the rectangles and the tori.

Gichev \cite[Theorem 3]{G} proved that on the $n$-$\dim$ unit sphere $\s^n$,
\begin{equation}\label{eq:Gichev}
H_2(\s^n)=\limsup_{\lambda\to\infty}\frac{\HH^{n-1}(\NN(u))}{\sqrt\lambda}=\Vol(\s^{n-1})=\frac{\Gamma((n+1)/2)}{\Gamma(n/2)\sqrt\pi}\Vol(\s^n),
\end{equation}
in which $\Vol(\s^{n-1})$ is the volume of $\s^{n-1}$. For example, $H_2(\s^2)=2\pi=\Vol(\s^2)/2$. In fact, the eigenfunctions on the sphere are spherical harmonics (i.e. homogeneous harmonic polynomials restricted to the sphere) and Gichev proved a stronger result that 
$$\HH^{n-1}(\NN(u))\le k\Vol(\s^{n-1}),$$
in which $k$ is the homogeneous degree of $u$. Moreover, the equation in the above inequality is obtained by the Gaussian beams (i.e. highest weight spherical harmonics). One then deduces \eqref{eq:Gichev} by observing that $\lambda=k(k+n-1)$. In the same paper \cite[Page 563]{G}, Gichev conjectured that on $\s^2$,
$$H_1(\s^2)=\liminf_{\lambda\to\infty}\frac{\HH^1(\NN(u))}{\sqrt\lambda}=4=\frac1\pi\Area(\s^2),$$
and the limit is achieved by the zonal harmonics.

Our main result is to provide the case in the disc for which \textit{both} of the sharp constants $H_1$ and $H_2$ in \eqref{eq:infsup} are explicitly proved, see Theorem \ref{thm:disc}; moreover, we observe the geometric properties of the eigenfunctions which achieve the $\limsup\HH^1(\NN(u))/\sqrt\lambda$ as $\lambda\to\infty$, see Corollary \ref{cor:geod}.

Let $\B_1=\{(x,y)\in\R^2:x^2+y^2<1\}$ be the unit disc. Consider the eigenfunctions with Dirichlet boundary condition. In polar coordinates $\{(r,\theta):0\le r\le1,0\le\theta<2\pi\}$, the real-valued Dirichlet eigenfunctions are
$$u_{k,s}(r,\theta)=J_k\left(\sqrt{\lambda_{k,s}}\cdot r\right)\sin(k\theta+\theta_0),\quad\text{where }k=0,1,2,3,...\text{ and }s=1,2,3,...$$
Here, $J_k$ is the $k$-th Bessel function and $\lambda_{k,s}=j_{k,s}^2$, where $j_{k,s}$ is the $s$-th nonnegative zero of $J_k$, and $\theta_0\in[0,2\pi)$. So all the Dirichlet eigenvalues are the squares of zeros of Bessel functions. In particular, the eigenvalues $j^2_{k,s}$ are distinct for different values of $k$ and $s$, (c.f. \cite[Section 15.23]{W}) and the multiplicity of $j^2_{k,s}$ is two with eigenspace spanned by 
$$J_k\left(\sqrt{\lambda_{k,s}}\cdot r\right)\sin(k\theta)\quad\text{and}\quad J_k\left(\sqrt{\lambda_{k,s}}\cdot r\right)\cos(k\theta).$$
The nodal set (in the interior of the disc) of the eigenfunction with eigenvalue $\lambda_{k,s}=j_{k,s}^2$ is a collection of $2k$ radials (i.e. $k$ diameters) and $s-1$ concentric circles with radii $j_{k,l}/j_{k,s}$, $l=1,...,s-1$.  Our main theorem states that
\begin{thm}\label{thm:disc}
In the disc $\B_1$, let $u$ be a Dirichlet eigenfunction with eigenvalue $\lambda$. Then
$$H_1(\B_1)=\liminf_{\lambda\to\infty}\frac{\HH^1(\NN(u))}{\sqrt\lambda}=1,$$
in which the limit is achieved by the eigenfunctions $u_{0,s}$ as $s\to\infty$; and
$$H_2(\B_1)=\limsup_{\lambda\to\infty}\frac{\HH^1(\NN(u))}{\sqrt\lambda}=2,$$
in which the limit is achieved by the eigenfunctions $u_{k,1}$ as $k\to\infty$.
\end{thm}

\begin{rmk}
The same results as in Theorem \ref{thm:disc} hold for Neumann eigenfunctions in the disc. See the discussion in \S\ref{sec:Neu}.
\end{rmk}

A simple dilation gives
\begin{cor}\label{cor:discR}
In the disc $\B_R=\{(x,y)\in\R^2:x^2+y^2<R^2\}$, let $u$ be a Dirichlet eigenfunction with eigenvalue $\lambda$. Then
$$H_1(\B_R)=\liminf_{\lambda\to\infty}\frac{\HH^1(\NN(u))}{\sqrt\lambda}=R^2=\frac{1}{\pi}\,\Area(\B_R)$$
and
$$H_2(\B_R)=\limsup_{\lambda\to\infty}\frac{\HH^1(\NN(u))}{\sqrt\lambda}=2R^2=\frac{2}{\pi}\,\Area(\B_R).$$
\end{cor}

\begin{rmk}
Some bounds of $c_1$ in \eqref{eq:Yau} on Riemannian surfaces are previously known, with which one can have some non-sharp estimates of $H_1$. In a Euclidean domain $\Omega\subset\R^2$, Br\"uning-Gromes \cite[Equation (4)]{BG} proved that
$$H_1(\Omega)\ge\frac{1}{2j_{0,1}}\,\Area(\Omega),$$
where, as before, $j_{0,1}\approx2.4048$ is the first zero of the Bessel function $J_0$. On a smooth Riemannian surface $\M$, Savo \cite[Theorem 13]{S} proved that 
$$H_1(\M)\ge\frac{1}{11}\,\Area(\M).$$
So our calculation of $H_1$ in Corollary \ref{cor:discR} can be regarded as the sharp improvement of these results applied to the disc.
\end{rmk}

The other problem in question is to characterize the (possible) geometric properties of the eigenfunctions that achieve  $\liminf$ or $\limsup$ of $\HH^1(\NN(u))/\sqrt\lambda$ as $\lambda\to\infty$. Here, we make the observation that the nodal set of 
$$u_{k,1}=J_k\left(\sqrt{\lambda_{k,1}}\cdot r\right)\sin(k\theta+\theta_0),$$ 
is a collection of $k$ diameters that pass through the origin. Hence,
\begin{cor}\label{cor:geod}
In the disc $\B_R$, 
$$H_2(\B_R)=\limsup_{\lambda\to\infty}\frac{\HH^1(\NN(u))}{\sqrt\lambda}$$
is saturated by a sequence of eigenfunctions whose nodal curves in the interior are geodesics, i.e. pieces of straight lines.
\end{cor}

Recall that on $\s^n$ proved by Gichev \cite{G}, $\limsup\HH^{n-1}(\NN(u))/\sqrt\lambda$ is saturated by the Gaussian beams, whose nodal sets are totally geodesic. Based on these evidence, we propose the following problem.
\begin{pblm}\label{pblm}
In what manifold $(\M,g)$, one has that
$$H_2(\M)=\limsup_{\lambda\to\infty}\frac{\HH^{n-1}(\NN(u))}{\sqrt\lambda}=\lim_{k\to\infty}\frac{\HH^{n-1}(\NN(u_k))}{\sqrt{\lambda_k}}$$
for a sequence of eigenfunctions $\{u_k\}_{k=1}^\infty$ with eigenvalues $\lambda_k$ such that the nodal sets of $u_k$ are totally geodesic in the interior of $\M$?
\end{pblm}

\begin{rmk}
The answer to Problem \ref{pblm} is positive on the spheres by Gichev \cite{G} and in the disc by Corollary \ref{cor:geod}. In the irrational rectangles (see \S\ref{sec:rectangles}), the nodal curves of all eigenfunctions are geodesics so the answer to Problem \ref{pblm} is trivially positive. It would be interesting to see if Problem \ref{pblm} holds in other rectangles (or tori). On a general manifold, the answer to Problem \ref{pblm} is not known and in fact it is not even known whether there exists a sequence of eigenfunctions whose nodal sets are totally geodesic. 
\end{rmk}

In all the manifolds that we consider in this paper (irrational rectangles and tori, spheres, and discs), we have that $H_1(\M)<H_2(\M)$. So a natural question follows as

\begin{pblm}\label{pblm:H}
In what manifold $(\M,g)$, $H_1(\M)=H_2(\M)$?
\end{pblm}

\begin{rmk}
In the case when $H_1(\M)=H_2(\M)$, there is a unique limit of $\HH^{n-1}(\NN(u))/\sqrt{\lambda}$ as $\lambda\to\infty$ for all the eigenfunctions. This is not known to be positive on any manifold with dimension higher than one. 

On an analytic manifold, one can extend the Laplacian eigenfunctions to a complex neighborhood of the manifold. In \cite[Corollary 1.2]{Z1}, Zelditch showed that on an analytic manifold with ergodic geodesic flow, there is a full density subsequence of eigenfunctions for which $\left[\NN(u^\C)\right]/\sqrt{\lambda}$
has a unique limit as $\lambda\to\infty$. Here, $u^\C$ denotes the complex extension of the eigenfunctions $u$ and $\left[\NN(u^\C)\right]$ denotes the complex hypersurface measure of the nodal set of $u^\C$. Even though this result is for the complex extensions of a full density subsequence of eigenfunctions, it suggests that on manifolds with ergodic geodesic flow (e.g. negatively curved manifolds), the answer to Problem \ref{pblm:H} might be positive.
\end{rmk}

\section{Proof of Theorem \ref{thm:disc}}
\subsection{The irrational rectangles and tori}\label{sec:rectangles}
Before proving Theorem \ref{thm:disc}, we discuss the proof of \eqref{eq:BG} in an irrational rectangle $R=\{(x,y)\in\R^2:0\le x\le a,0\le y\le b\}$, where $a^2/b^2$ is irrational. (This is observed in \cite{BG}.) We then make some remarks about finding $H_1$ and $H_2$ in more general rectangles and tori. 

The Dirichlet eigenfunctions in $R$ have the form
\begin{equation}\label{eq:toral}
u_{k,j}(x,y)=\sin\left(\frac{\pi k}{a}x\right)\sin\left(\frac{\pi j}{b}y\right)
\end{equation}
with the eigenvalues
$$\lambda_{k,j}=\left(\frac{\pi k}{a}\right)^2+\left(\frac{\pi j}{b}\right)^2,\quad\text{where }k,j=1,2,3...$$
Notice that if $a^2/b^2$ is irrational, then all the eigenvalues are simple. Indeed, if $\lambda_{\tilde k,\tilde j}=\lambda_{k,j}$ for another pair $(\tilde k,\tilde j)$, $\tilde k,\tilde j=1,2,3...$, then 
$$\frac{k^2-\tilde k^2}{a^2}=\frac{\tilde j^2-j^2}{b^2},$$
which forces $(\tilde k,\tilde j)=(k,j)$ since $a^2/b^2$ is irrational. The nodal set $\NN(u_{k,j})$ (in the interior of $R$) consists of $(k-1)$ line segments of length $b$ and $(j-1)$ line segments of length $a$. So the nodal length of $u_{k,j}$ is
$$\HH^1(\NN(u_{k,j}))=(k-1)b+(j-1)a.$$
Hence,
$$\frac{\HH^1(\NN(u_{k,j}))}{\sqrt{\lambda_{k,j}}}=\frac{\Area(R)}{\pi}\times\frac{(k-1)b+(j-1)a}{\sqrt{(kb)^2+(ja)^2}}.$$
One then sees from $\sqrt{p^2+q^2}\le p+q\le\sqrt2\sqrt{p^2+q^2}$ for $p,q\ge0$ that 
\begin{equation}\label{eq:H1R}
H_1(R)=\liminf_{\lambda\to\infty}\frac{\HH^1(\NN(u_{k,j}))}{\sqrt{\lambda_{k,j}}}=\frac1\pi\Area(R),
\end{equation}
in which the limit is achieved by $u_{1,j}$ as $j\to\infty$ and by $u_{k,1}$ as $k\to\infty$, and
\begin{equation}\label{eq:H2R}
H_2(R)=\limsup_{\lambda\to\infty}\frac{\HH^1(\NN(u_{k,j}))}{\sqrt{\lambda_{k,j}}}=\frac{\sqrt2}{\pi}\Area(R),
\end{equation}
in which the limit is achieved by $u_{k,j}$ such that $k/j\to a/b$ as $k,j\to\infty$.

\begin{rmk}\hfill
\begin{itemize}
\item The above proof holds with little modification for Neumann eigenfunctions in these irrational rectangles.
\item In a more general rectangles, the eigenvalues may not be simple and can have high multiplicity, e.g. in the rectangle $[0,\pi]\times[0,\pi]$, there are eigenvalues $\lambda$ with multiplicity of the order $\sqrt{\log\lambda}$ as $\lambda\to\infty$. In the case of high multiplicity, one has to estimate the precise nodal lengths of linear combinations of eigenfunctions of the form \eqref{eq:toral}. These linear combinations have complex nodal portraits and the problem of finding $H_1$ and $H_2$ becomes challenging.
\item In the disc, the eigenvalues have multiplicity two. Therefore, in the following subsection, we can use explicit formulae to deduce $H_1$ and $H_2$.
\end{itemize}
\end{rmk}

If we identify the two opposing sides of the rectangle $R=[0,a]\times[0,b]$ and define the torus $\T$ (i.e. without boundary), then the real-valued eigenfunctions are spanned by
$$\sin\left(\frac{2\pi k}{a}x\pm\frac{2\pi j}{b}y\right)\quad\text{and}\quad\cos\left(\frac{2\pi k}{a}x\pm\frac{2\pi j}{b}y\right)$$
with eigenvalue
$$\left(\frac{2\pi k}{a}\right)^2+\left(\frac{2\pi j}{b}\right)^2,\quad\text{where }k,j=0,1,2,3...$$
Given that $a^2/b^2$ is irrational, the eigenvalues are not simple (except when $k=j=0$) but their multiplicity is uniformly bounded by $4$. A similar argument as in the corresponding irrational rectangle shows that the same results of $H_1$ and $H_2$ in \eqref{eq:H1R} and \eqref{eq:H2R} hold on the torus. However, $H_1$ and $H_2$ remain unknown on other tori, for the same reason as described in the above remark.

\subsection{Proof of Theorem \ref{thm:disc}}\label{sec:disc}
Now we prove Theorem \ref{thm:disc}. Recall that the Dirichlet eigenfunction with eigenvalue $\lambda_{k,s}=j_{k,s}^2$ has the form
$$u_{k,s}=J_k\left(\sqrt{\lambda_{k,s}}\cdot r\right)\sin(k\theta+\theta_0),\quad\text{where }k=0,1,2,3,...\text{ and }s=1,2,3,...$$
Since the nodal length $\HH^1(\NN(u_{k,s}))$ is independent of $\theta_0$ here, we assume $\theta_0=0$ without loss of generality. The nodal set of $u_{k,s}$ is a collection of $k$ diameters and $s-1$ concentric circles with radii $j_{k,l}/j_{k,s}$, $l=1,...,s-1$. In particular, $\NN(u_{0,s})$ consists of circles only and $\NN(u_{k,1})$ consists of diameters only. 

Here, we provide the graphs of the nodal curves of some eigenfunctions with different nodal portrait.
\begin{figure}[h!]
\includegraphics[width=4cm]{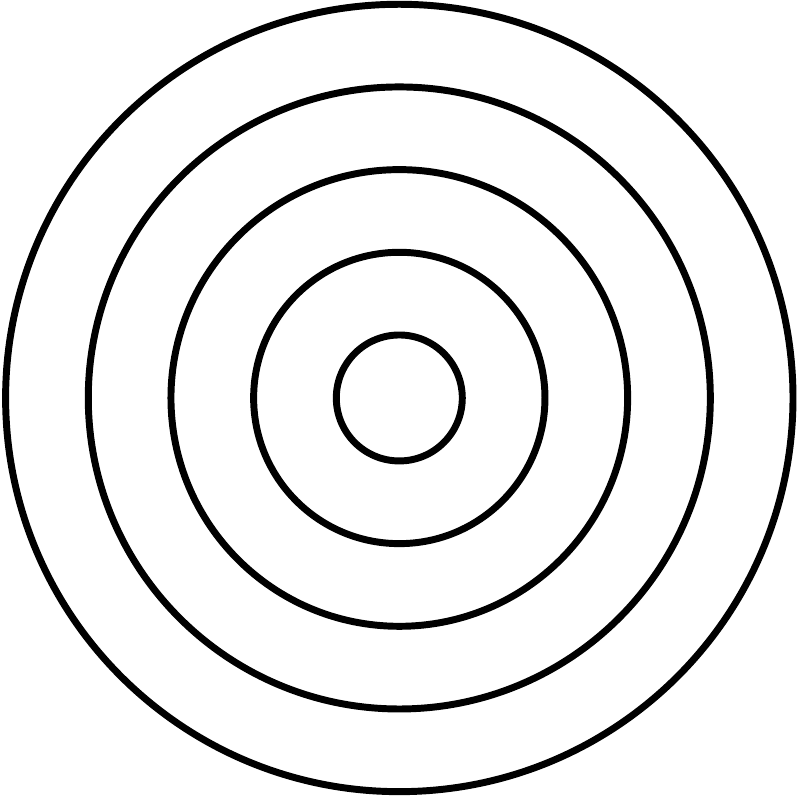}\quad\quad\quad\quad
\includegraphics[width=4cm]{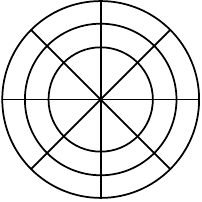}\quad\quad\quad\quad
\includegraphics[width=4cm]{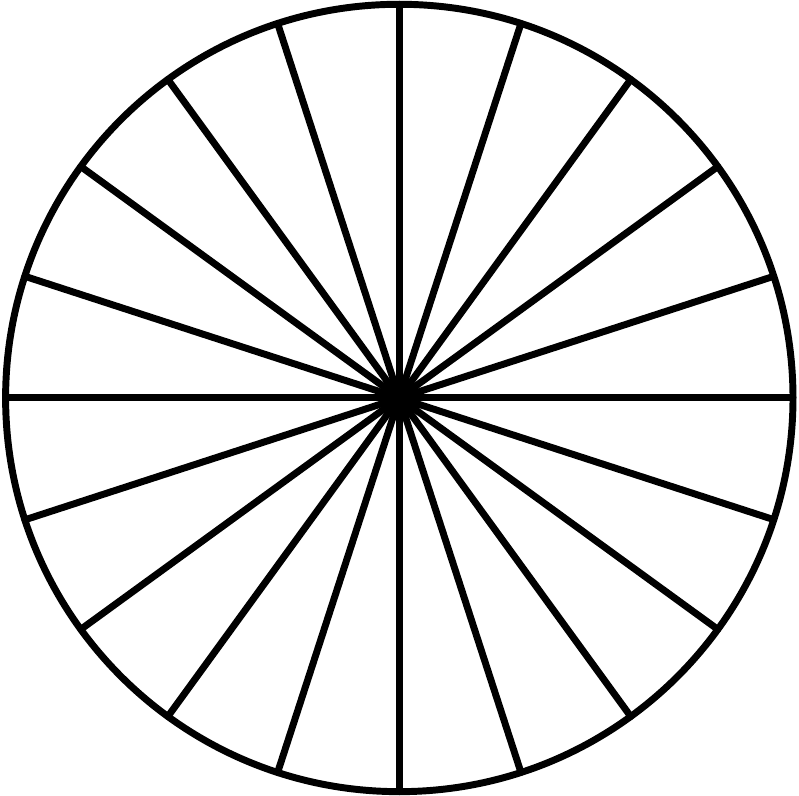}
\end{figure}

From left to right: $u_{0,5}$, $u_{4,3}$, and $u_{10,1}$. Their eigenvalues are approximately $14^2$ and one can check that $\HH^1(\NN(u_{0,5}))<\HH^1(\NN(u_{4,3}))<\HH^1(\NN(u_{10,1}))$, which is a reflection of Theorem \ref{thm:disc}.

To estimate 
$$H_1(\B_1)=\liminf_{\lambda\to\infty}\frac{\HH^1(\NN(u))}{\sqrt\lambda}\quad\text{and}\quad\quad H_2(\B_2)=\limsup_{\lambda\to\infty}\frac{\HH^1(\NN(u))}{\sqrt\lambda},$$
we pick any subsequence of $\{u_{k,s}\}$ such that
$$\lim_{\lambda\to\infty}\frac{\HH^1(\NN(u))}{\sqrt\lambda}\quad\text{exists},$$
and divide into three cases.

\begin{itemize}
\item Case 1: $k$ tends to infinity and $s$ is bounded;
\item Case 2: $s$ tends to infinity and $k$ is bounded;
\item Case 3: $k$ and $s$ both tend to infinity.
\end{itemize}

\textbf{Case 1}. As $\lambda=j_{k,s}^2\to\infty$, $s$ is bounded so $k\to\infty$. First set $s=1$. Then the nodal set of $u_{k,1}$ is the union of $2k$ radials, that is,
$$\HH^1(\NN(u_{k,1}))=2k,$$
and
$$\lim_{k\to\infty}\frac{\HH^1(\NN(u_{k,1}))}{\sqrt{\lambda_{k,1}}}=\lim_{k\to\infty}\frac{2k}{j_{k,1}}=2.$$
Here, we use the fact that from \cite[Equation 9.5.14, pp. 371]{AS},
\begin{equation}\label{eq:kinfty}
j_{k,1}=k+O\left(k^\frac13\right)\quad\text{as }k\to\infty.
\end{equation}
This argument works for all the subsequences of $\{u_{k,s}\}$ for which $s$ is bounded and $k\to\infty$. Indeed, If $s$ is bounded by $M$, then
$$2k\le\HH^1(\NN(u_{k,s}))=2\pi\sum_{l=1}^{s-1}\frac{j_{k,l}}{j_{k,s}}+2k\le2\pi M+2k.$$
So by squeezing,
$$\lim_{k\to\infty,\ s\text{ bounded}}\frac{\HH^1(\NN(u_{k,s}))}{\sqrt{\lambda_{k,s}}}=2.$$
Here, we need to use a formula \cite[Equation 9.5.22, pp. 371]{AS} that if $s$ is bounded, then
$$j_{k,s}=k+o(k)\quad\text{as }k\to\infty.$$ 

\textbf{Case 2}. As $\lambda=j_{k,s}^2\to\infty$, $k$ is bounded so $s\to\infty$. First set $k=0$. Then the nodal set of $u_{0,s}$ is the union of $s-1$ concentric circles with radii $j_{0,l}/j_{0,s}$, $l=1,...,s-1$, that is,
$$\HH^1(\NN(u_{0,s}))=2\pi\sum_{l=1}^{s-1}\frac{j_{0,l}}{j_{0,s}},$$
and
$$\lim_{s\to\infty}\frac{\HH^1(\NN(u_{0,s}))}{\sqrt{\lambda_{0,s}}}=2\pi\lim_{s\to\infty}\sum_{l=1}^{s-1}\frac{j_{0,l}}{j_{0,s}^2}=2\pi\lim_{s\to\infty}\frac{\sum_{l=1}^s\left(l-\frac14\right)\pi}{\left[\left(s-\frac14\right)\pi\right]^2}=1.$$
Here, we use the fact that from \cite[Equation 9.5.12, pp. 371]{AS}: If $k\ll s$, then
\begin{equation}\label{eq:sinfty}
j_{k,s}=\left(s+\frac k2-\frac14\right)\pi+O\left(s^{-1}\right)\quad\text{as }s\to\infty.
\end{equation}
This argument works for all the subsequences of $\{u_{k,s}\}$ for which $k$ is bounded and $s\to\infty$. Indeed, If $k$ is bounded by $M$, then $\NN(u_{k,s})$ contains $s-1$ circles and at most $M$ diameters. Hence,
$$2\pi\sum_{l=1}^{s-1}\frac{j_{k,l}}{j_{k,s}}\le\HH^1(\NN(u_{k,s}))=2\pi\sum_{l=1}^{s-1}\frac{j_{k,l}}{j_{k,s}}+2k\le2\pi\sum_{l=1}^s\frac{j_{k,l}}{j_{k,s}}+2M.$$
So by squeezing, using \eqref{eq:sinfty} again, we have that
$$\lim_{s\to\infty,\ k\text{ bounded}}\frac{\HH^1(\NN(u_{k,s}))}{\sqrt{\lambda_{k,s}}}=1.$$

\textbf{Case 3}. As $\lambda=j_{k,s}^2\to\infty$, $k,s\to\infty$. Suppose that in such a subsequence
$$\lim_{\lambda\to\infty}\frac{\HH^1(\NN(u_{k,s}))}{\sqrt{\lambda_{k,s}}}=p.$$
Our goal is then to prove that $1\le p\le2$.

We need the uniform bounds of the zeros $j_{k,s}$ of Bessel functions as $k,s\to\infty$. By \cite[Equation (1) in Theorem 1]{Bre}, we have that
$$j_{k,s}>k+\frac23|a_{s-1}|^\frac32\quad\text{for }k=0,1,2,3,...\text{ and }s=1,2,3,...$$
Here, $a_s$ is the $s$-th negative zero of the Airy function. By \cite[Equations 10.4.94 and 10.4.105]{AS}, we have that as $s\to\infty$,
$$a_s=-\left[\frac{3\pi(4s-1)}{8}\right]^\frac23\left[1+O\left(s^{-2}\right)\right].$$
Hence,
\begin{equation}\label{eq:jlowerbd}
j_{k,s}>k+\frac23|a_{s-1}|^\frac32=k+\pi s+O\left(s^{-1}\right).
\end{equation}
We now proceed to prove the upper bound that $p\le2$. Using \eqref{eq:jlowerbd},
$$\frac{\HH^1(\NN(u_{k,s}))}{\sqrt{\lambda_{k,s}}}=\frac{2\pi\sum_{l=1}^{s-1}\frac{j_{k,l}}{j_{k,s}}+2k}{j_{k,s}}<\frac{2\pi s+2k}{k+\pi s+O\left(s^{-1}\right)}\le2\quad\text{as }k,s\to\infty.$$
We then prove the lower bound that $p\ge1$. By \cite[Equation (2) in Theorem 1]{Bre}, we have that
\begin{equation}\label{eq:jupperbd}
j_{k,s}<\frac\pi2k+\frac23|a_s|^\frac32=\frac\pi2k+\pi s+O\left(s^{-1}\right)\quad\text{for }k=0,1,2,3,...\text{ and }s=1,2,3,...
\end{equation}
Using \eqref{eq:jlowerbd} and \eqref{eq:jupperbd}, we compute that
\begin{eqnarray*}
\sum_{l=1}^{s-1}\frac{j_{k,l}}{j_{k,s}}&\ge&\frac{1}{\frac\pi2k+\pi s+O\left(s^{-1}\right)}\sum_{l=1}^{s-1}\left[k+\pi l+O\left(l^{-1}\right)\right]\\
&\ge&\frac{\frac\pi2(s-1)s+(s-1)k+O(s)}{\frac\pi2k+\pi s+O\left(s^{-1}\right)}\\
&\ge&\frac s2,\quad\text{if $s$ and $k$ are large enough}.
\end{eqnarray*}

\begin{rmk}
Notice that $j_{k,l}/j_{k,s}$, $l=1,...,s-1$, are fractions which are distributed in the interval $[0,1]$. If $k\ll s$, then by the asymptotic formula \eqref{eq:sinfty}, these fractions are rather equidistributed. So the above inequality is natural in this case. If $s\ll k$, then by \eqref{eq:kinfty} and \eqref{eq:jupperbd}, one sees that $j_{k,1}/j_{k,s}\gtrsim k/(\pi k/2+\pi s)>1/2$ and therefore $j_{k,l}/j_{k,s}>1/2$ for all $l=1,...,s-1$. So the above inequality is natural in this case as well. The above inequality in fact shows that it is true for all sufficiently large $k$ and $s$.
\end{rmk}

Now by \eqref{eq:jupperbd} again,
$$\frac{\HH^1(\NN(u_{k,s}))}{\sqrt{\lambda_{k,s}}}=\frac{2\pi\sum_{l=1}^{s-1}\frac{j_{k,l}}{j_{k,s}}+2k}{j_{k,s}}\ge\frac{\pi s+2k}{\frac\pi2k+\pi s+O\left(s^{-1}\right)}\ge1\quad\text{as }k,s\to\infty.$$
Hence, the lower bound is proved.

\subsection{Neumann eigenfunctions}\label{sec:Neu}
The Neumann eigenfunctions in $\B_1$ can be written as
$$v_{k,s}(r,\theta)=J_k\left(\sqrt{\mu_{k,s}}\cdot r\right)\sin(k\theta+\theta_0),\quad\text{where }k=0,1,2,3,...\text{ and }s=1,2,3,...$$
Here, $J_k$ is the $k$-th Bessel function and $\mu_{k,s}=(j'_{k,s})^2$, where $j'_{k,s}$ is the $s$-th nonnegative zero of $J_k'$. So all the Neumann eigenvalues are the squares of zeros of the derivatives of Bessel functions. 

Here for Neumann eigenfunctions, we could repeat the argument in the previous subsection, using instead the estimates of $j_{k,s}'$. However, notice that Neumann eigenfunctions in $\B_1$ extends to $\R^2$ and in fact defines a Dirichlet eigenfunction in a slightly larger disc. So we can estimate the nodal set of Neumann eigenfunctions by Corollary \ref{cor:discR} for Dirichlet eigenfunctions.

Indeed, the zeros $j'_{k,s}$ and $j_{k,s}$ interlace according to
\begin{equation}\label{eq:interlace}
k\le j_{k,s}'<j_{k,s}<j_{k,s+1}'.
\end{equation}
See \cite[Equation 9.5.2]{AS}. Using this relation, we see that $v_{k,s}$ extends from $\B_1$ to $\B_R$ as a Dirichlet eigenfunction $u_{k,s}$ with
\begin{equation}\label{eq:R}
R=\frac{j_{k,s}}{j'_{k,s}}\to1\quad\text{as }k\text{ or }s\to\infty.
\end{equation}
Now the nodal set of $v_{k,s}$ in $\B_1$ and the nodal set of $u_{k,s}$ in $\B_R$ differ by the $2k$ radials in $\B_R\setminus\B_1$. That is,
$$\HH^1(\NN(v_{k,s}))=\HH^1(\NN(u_{k,s}))-2k(R-1).$$
Hence,
$$\frac{\HH^1(\NN(v_{k,s}))}{j'_{k,s}}=\frac{\HH^1(\NN(u_{k,s}))-2k(R-1)}{j_{k,s}}\cdot R.$$
By \eqref{eq:interlace} and \eqref{eq:R}, we see that
$$\frac{2k(R-1)}{j_{k,s}}\to0\quad\text{as }j_{k,s}\to\infty.$$
Then applying Corollary \ref{cor:discR} and again \eqref{eq:R}, we have that
$$\liminf_{\lambda\to\infty}\frac{\HH^1(\NN(v))}{\sqrt\mu}=1\quad\text{and}\quad\limsup_{\lambda\to\infty}\frac{\HH^1(\NN(v))}{\sqrt\mu}=2.$$

\section*{Acknowledgments}
XH wants to thank Stephen Breen, Andrew Hassell, Hamid Hezari, and Steve Zelditch for all the discussions that are related to this article, in particular, Problem \ref{pblm}; XH also wants to thank Ze\'ev Rudnick for informing him the results in Gichev \cite{G} and Werner Horn for his translation of Br\"uning-Gromes \cite{BG}.

\end{document}